\documentclass[letterpaper, 10 pt, conference]{ieeeconf}

\IEEEoverridecommandlockouts                              % This command is only
 \pdfoutput=1
\usepackage{thmtools,thm-restate}
\usepackage{amsmath, amssymb}
\usepackage{amsfonts}
\usepackage{graphicx}
\usepackage{enumerate}
\usepackage{ifthen}
\usepackage{multicol}
\usepackage{float}
\usepackage{cite}
\usepackage{fancyhdr}
\usepackage[usenames, dvipsnames]{color}
\graphicspath{{figures/}}
\usepackage{mathtools}
\usepackage{multicol}
\usepackage{xcolor}

\usepackage{svg}
 \usepackage{cancel} %To be able to cross out text
 \usepackage[normalem]{ulem}

\newtheorem{theorem}{Theorem}
\newtheorem{lemma}[theorem]{Lemma}
\newtheorem{proposition}[theorem]{Proposition}
\newtheorem{corollary}[theorem]{Corollary}
\newtheorem{rem}{Remark}
\newtheorem{example}{Example}
\newtheorem{ass}{Assumption}
\newtheorem{definition}{Definition}

\newboolean{showcomments}
\setboolean{showcomments}{true}

\newcommand{\bmat}[1]{\begin{bmatrix}
#1
\end{bmatrix}}

\newcommand{\todo}[1]{  \ifthenelse{\boolean{showcomments}}
{\textcolor{ForestGreen}{TO DO:  #1}}{}}
\newcommand{\suggest}[1]{\ifthenelse{\boolean{showcomments}}
{\textcolor{Orange}{(Suggestion: #1)}}{}}
\newcommand{\alain}[1]{\ifthenelse{\boolean{showcomments}}
{\textcolor{Blue}{(Alain says: #1)}}{}}
\newcommand{\jonas}[1]{\ifthenelse{\boolean{showcomments}}
{\textcolor{ForestGreen}{(Jonas says: #1)}}{}}
\newcommand{\kristian}[1]{\ifthenelse{\boolean{showcomments}}
{\textcolor{Blue}{(Kristian says: #1)}}{}}
\newcommand{\emma}[1]{\ifthenelse{\boolean{showcomments}}
{\textcolor{VioletRed}{(Emma says: #1)}}{}}
\newcommand{\ifneeded}[1]{\ifthenelse{\boolean{showcomments}}
{\textcolor{Gray}{#1}}{}}

\newboolean{showedit}
\setboolean{showedit}{true}
\newcommand{\edit}[1]{\ifthenelse{\boolean{showedit}}
{\textcolor{Blue}{#1}}{}}
\newcommand{\draft}[1]{\ifthenelse{\boolean{showedit}}
{\textcolor{gray}{#1}}{}}

\let\labelindent\relax
\usepackage[shortlabels]{enumitem}

\usepackage{array}
\newcolumntype{L}[1]{>{\raggedright\let\newline\\\arraybackslash\hspace{0pt}}m{#1}}
\newcolumntype{C}[1]{>{\centering\let\newline\\\arraybackslash\hspace{0pt}}m{#1}}
\newcolumntype{R}[1]{>{\raggedleft\let\newline\\\arraybackslash\hspace{0pt}}m{#1}}

\usepackage{physics} %to enable use of \norm

\usepackage[font = small]{caption}
\usepackage{subcaption}

\usepackage{tikz}
 \usetikzlibrary{plotmarks}
\usepackage{tikz}
\usepackage{pgfplots}
\pgfplotsset{compat=newest}
\usetikzlibrary{patterns}
\usetikzlibrary{decorations.text}
\usepgfplotslibrary{fillbetween}

\newcommand{\ts}{\textsuperscript}

\providecommand{\figref}{}
\renewcommand{\figref}[1]{Fig.~\ref{#1}}
\providecommand{\secref}{}
\renewcommand{\secref}[1]{Sec.~\ref{#1}}
\providecommand{\propref}{}
\renewcommand{\propref}[1]{\textit{Proposition}~\ref{#1}}

\usepackage{lipsum}
\usepackage{mathtools}

\title{\LARGE \bf A closed-loop design for scalable high-order consensus}

\author{{Jonas Hansson and Emma Tegling} 
 \thanks{The authors are with the Department of Automatic Control and the ELLIIT Strategic Research Area at 
        Lund University, Lund, Sweden. Email: \{{\tt\small{jonas.hansson, emma.tegling}\}@control.lth.se}}\thanks{This work was partially funded by Wallenberg AI, Autonomous Systems and Software Program (WASP) funded by the Knut and Alice Wallenberg Foundation and the Swedish Research Council through Grant 2019-00691. }}

\begin{document}
\maketitle
\begin{abstract}
    This paper studies the problem of coordinating a group of $n$\ts{th}-order integrator systems. As for the well-studied conventional consensus problem, we consider 
    linear and distributed control with only local and relative measurements.  
     We propose a closed-loop dynamic that we call  \emph{serial consensus} and prove it achieves $n\ts{th}$ order consensus regardless of model order and underlying network graph. This alleviates {an important }scalability limitation in conventional consensus dynamics of order $n\ge 2$, whereby they may lose stability if the underlying network grows.  The  
     distributed control law which achieves the desired closed loop dynamics is shown to be localized and obey the limitation to relative state measurements.  Furthermore, through use of the small-gain theorem, the serial consensus system is shown to be robust to both model and feedback uncertainties.   We illustrate the theoretical results through examples.   
\end{abstract}

\section{Introduction}
Properties of dynamical systems over networks have been a subject of significant research over the last two decades. A problem of interest the coordination of agents in a network through localized feedback, leading to the prototypical  distributed consensus dynamics, studied early on by~\cite{FaxMurray,OlfatiSaber2004,Jadbabaie2003}. Over the years, it has become clear that the structural constrains imposed by the network topology in consensus problems often lead to fundamentally poor dynamic behaviors in large networks. This concerns controllability~\cite{Pasqualetti2014}, performance~\cite{Bamieh2012,SiamiMotee2015} and disturbance propagation~\cite{swaroop1996stringstability,seiler2004disturbancep_propagation}, but, as recently highlighted in~\cite{Tegling2022_Scalefrag}, also stability. 
The poor stability properties characterized in earlier work~\cite{Tegling2022_Scalefrag} (which motivate the present work) apply to higher-order consensus, where the local dynamics of each agent is modeled
as an $n$th order integrator, with $n\geq2$, and the control is a weighted average of neighbors' relative states. 
This is a theoretical generalization of first-order consensus~\cite{Jiang2009}, but is also relevant in practice. For example, a model where $n=3$ and thus has consensus in position, velocity and acceleration, can capture flocking behaviors~\cite{Ren2006}. 

More specifically, \cite{Tegling2022_Scalefrag} shows that conventional high-order consensus $(n\geq 3)$ is not \textit{scalably stable} for many growing graph structures. When the network grows beyond a certain size, stability is lost. The same holds for second-order consensus ($n = 2$) in, for example, directed ring graphs, as also described in~\cite{Studli2017}. 
To address this lack of scalable stability we propose an alternative generalization of the first-order consensus dynamics, which we prove achieves scalable stability for any model order~$n$. 

To illustrate our proposed controller, consider the conventional second-order consensus system where the controller 
$u(t)=-L_1\dot{x}(t)-L_2x(t)+u_\mathrm{ref}(t)$, with $L_{1,2}$ being weighted graph Laplacians, is used to achieve the closed loop
\begin{equation}
    \ddot{x}(t)=-L_2\dot{x}(t)-L_1x+u_\mathrm{ref}(t).
    \label{eq:secondorderconventional}
\end{equation}
While for first-order consensus ($\dot{x} = Lx+u_\mathrm{ref}(t) $), a sufficient condition for convergence to consensus is that the graph underlying the graph Laplacian $L$ contains a connected spanning tree. However, this no longer suffices when~{$n\geq2$} as in~\eqref{eq:secondorderconventional}. 
 Therefore, we instead propose the following controller $u(t)=-(L_1+L_2)\dot{x}(t)-L_1L_2x(t)+u_\mathrm{ref}(t)$. The reason for this choice of controller is best illustrated by considering the resulting closed loop in the Laplace domain:
 \begin{equation}
     (sI+L_1)(sI+L_2)X(s)=U_\mathrm{ref}.
     \label{eq:secondorderserial}
 \end{equation}
For this system, like for the first-order case, it is sufficient that the graphs underlying $L_1$ and $L_2$ contain a connected spanning tree for the system to eventually coordinate in both~$x$ and its derivative~$\dot{x}$ (regardless of network size!). This closed loop system, which we will call \emph{serial consensus}, thus mimics one core property of the standard %regular 
consensus protocol, and can also be generalized to any order $n$.

The main results of this paper are proofs of some key properties of the proposed $n$\ts{th}-order serial consensus. The controller 
is proven to remain localized (within an $n$-hop neighborhood) and implementable through relative measurements. We also prove that the closed loop will achieve consensus in all $n$ states.
Furthermore, we study the robustness of the proposed closed loop and show that the system will still coordinate when subject to unstructured uncertainty.  The beneficial properties of the form~\eqref{eq:secondorderserial} (generalized to any order $n$) are thus not contingent on an idealized implementation.

The remainder of this paper is organized as follows.
We first introduce the $n$\ts{th} order consensus model and define our choice of control structure. Then the serial consensus system is defined and motivated. In \secref{sec:Results} we provide proofs for the stability and robustness of the serial consensus system. Our main results are then illustrated through examples in \secref{sec:Examples}. Lastly, we provide our Conclusions in \secref{Conclusion}.

%%%%%%%%%%%%%%%%%%%%%%%%%%%%%%%%
\section{Problem Setup}
We start by introducing some graph theory  before %and then commence with 
introducing the general $n$\ts{th} order consensus problem for which we propose the new serial consensus setup. We discuss its properties and then end with some useful definitions. 

\subsection{Network model and definitions}
Let $\mathcal{G}=\{\mathcal{V},\mathcal{E}\}$ denote a graph of size $N=|\mathcal{V}|$. The set $\mathcal{E}\subset \mathcal{V} \times \mathcal{V}$ denotes the set of edges. The graph can be equivalently represented by the adjacency matrix $W\in \mathbb R ^{N\times N}$ where $w_{i,j}>0 \iff (j,i)\in \mathcal{E}$. The graph is called \emph{undirected} if $W=W^T$. The graph contains a \emph{connected spanning tree} if for some $i\in \mathcal{V}$ there is a path from $i$ to any other vertex $j\in \mathcal{V}$.

Associated with a weighted graph we have the weighted graph Laplacian $L$ defined as
\begin{equation*}
    [L]_{i,j}= \left\{\begin{matrix}-w_{i,j},& ~~\mathrm{ if }\; i\neq j \\ 
                        \sum_{k\neq i} w_{i,k},& ~~\mathrm{ if }\; i=j\end{matrix}
    \right.
\end{equation*}
Under the condition that that the graph generating the graph Laplacian contains a connected spanning tree, $L$ will have a simple and unique  eigenvalue at $0$ 
and the remaining eigenvalues will lie strictly in the right half plane (RHP). 

We will also consider networks with a growing number of nodes. With $\mathcal{G}_N$ we denote a graph in a family $\{\mathcal{G}_N\}$, where $N$ is the size of the growing network. 

We will denote the space of all proper, real rational, and stable transfer matrices $\mathcal{R}\mathcal{H}_\infty$ and denote the $\mathcal{H}_\infty$ norm as $\|\cdot\|_{\mathcal{H}_\infty}$ following the notation in \cite{zhou1998essentials}.

\subsection{$n\ts{th}$ order consensus}
Let the system be modeled as $N$ agents with identical $n\ts{th}$ order integrator dynamics, i.e.
\begin{equation}
   \frac{\mathrm{d}^n x_i(t)}{\mathrm{d}t^n}=u_i(t), %\text{ for } i \in \mathcal{V}.
   \label{eq:n_order_integrator}
\end{equation}
for all $i \in \mathcal{V}$.
We will use the convention $x_i^{(0)}(t)=x_i(t)$ and $x_i^{(k)}(t)=\frac{\mathrm{d^k}}{\mathrm{d}t^k}x_i(t)$ to denote time derivatives. When clear we may omit the time argument for brevity.

In this paper we will consider the problem of synchronizing the agents and thus achieve a state of consensus.
\begin{definition}[$n$\ts{th} order consensus]
The multi-agent system \eqref{eq:n_order_integrator} is said to achieve ($n$\ts{th} order) consensus if  
$\lim_{t\rightarrow\infty}|x_i^{(k)}(t)-x_j^{(k)}(t)|=0$, for all $i,j\in \mathcal{V}$ and $k\in \{0,1\dots,n-1\}.$
\end{definition}

\subsection{Control structure}
A linear state feedback controller of \eqref{eq:n_order_integrator} can be written as 
\begin{equation}
u(t)= u_\mathrm{ref}(t)-\sum_{k=0}^{n-1}A_kx^{(k)}(t).
\label{eq:local_control}
\end{equation}
Where $u_\mathrm{ref}(t)\in \mathbb R ^N$ is a feedforward term and $A_k\in \mathbb R^{N\times N}$ represent the feedback of the $k$\ts{th} derivative.
We will restrict this class of controllers in three ways. The controllers
\begin{enumerate}[i)]
    \item can only use \emph{relative} feedback;
    \item have a limited gain;
    \item and depend on the local neighborhood of each agent.
\end{enumerate}
The limitation to relative feedback translates 
to the condition $A_k\mathbf{1}_N=0$ for all $k$, while a limited gain can be encoded by demanding that $\|A_k\|_\infty\leq c$. To  capture the notion of 
locality, consider the 
adjacency matrix $W$ representing the communication and measurement structure, which we here assume to be the same. That is, if $W_{i,j}=1$, then agent $i$ can directly receive or measure the relative distance to agent~$j$. Next, consider the non-negative matrix $W^q$. This matrix has the property that $[W^q]_{i,j}\neq 0$ if and only if there is a path of length $q$ from agent $j$ to agent $i$. Thus, if we want the controller to only depend on information that is at most $q$ steps away from each agent the following implication should hold: $\sum_{k=0}^{q} W^k_{i,j}=0 \implies [A_k]_{i,j}= 0$. Putting all the conditions together gives us a family of controllers that we will consider in this paper:
\begin{definition}[$q$-step implementable relative feedback]\label{def:q_implementable}
    A relative feedback controller of the form \eqref{eq:local_control} is $q$-step implementable with respect to the adjacency matrix $W$ and gain $c>0$ if $A_k\in\mathcal{A}^q(W,c)$ for all $k$,  where
    \begin{equation*}
        \mathcal{A}^q(W,c)= \begin{Bmatrix}
        A& \left|\begin{matrix}
        \left[\sum_{k=0}^{q} W^k\right]_{i,j}=0 \implies A_{i,j}= 0,\\ 
        A \mathbf{1}_N=0 , ~\|A\|_\infty\leq c
        \end{matrix}
        \right.
        \end{Bmatrix}
    .\end{equation*}  
\end{definition}
\vspace{1mm}

The conventional controller for achieving $n$\ts{th} order consensus can be realized as \eqref{eq:local_control} where each $A_k$ is given by a graph Laplacian, e.g. $A_k=L_k\in \mathcal{A}^1(W,c)$. In many cases these are also assumed to be the same such that $L_k=p_kL$ for some graph Laplacian $L$ and constants $p_k>0$.
\subsection{A Novel Design: Serial Consensus} \label{sec:serialconsensus}
We propose the following controller of \eqref{eq:n_order_integrator}, expressed in the Laplace domain, to achieve $n$\ts{th} order consensus
\begin{equation}
    U(s)=U_\mathrm{ref}(s)+ \left(s^nI -\prod_{k=1}^n (sI+L_k)\right)X(s),
    \label{eq:serial_consensus_controller}
\end{equation}
where $L_k$ are graph Laplacians and $U_\mathrm{ref}$ is the transformed reference signal.  In this case, it is more instructive to consider the closed-loop dynamics, which take the following form: %The resulting closed loop dynamics are then
\begin{definition}[$n$\ts{th} order serial consensus system]\label{def:serial_consensus}
For all $k\in\{1,2\dots, n\}$,  let $L_k$ be a weighted and directed graph Laplacian. 
The $n\ts{th}$-order serial consensus system is then % defined as
    \begin{equation}
        \left(\prod_{k=1}^{n} (sI+L_k)\right)X(s)=U_\mathrm{ref}(s).
        \label{eq:prod_laplac_def}
    \end{equation}
\end{definition}
\vspace{2mm}

\noindent We call this form serial consensus because the same closed loop dynamics can also be achieved by interconnecting $n$ first-order consensus systems in a series.
The closed-loop dynamics in \eqref{eq:prod_laplac_def} can also be transformed to state-space form by introducing the alternative variables $\Xi_k$ with the corresponding states $\xi_k$. These relate to $X$ through $\Xi_1=X(s)$, $\Xi_k=(sI+L_{k-1})\Xi_{k-1}$ for $k\in\{2,\dots,n~-~1\}$, and $s\Xi_n=-L_n\Xi_n+U_\mathrm{ref}$. This leads to the following continuous-time state-space representation
\begin{equation}
\begin{bmatrix}\dot{\xi}_1\\
\dot{\xi}_2\\
\vdots\\
\dot{\xi}_{n-1}\\
\dot{\xi}_n
\end{bmatrix}
\!=\!
\underbrace{\begin{bmatrix}-L_1 & I & \\
       & -L_2& \ddots \\
      &  & \ddots &I \\
      &  & &-L_{n} 
\end{bmatrix}
}_A \!\!
\begin{bmatrix}\xi_1\\
\xi_2\\
\vdots\\
\xi_{n-1}\\
\xi_n
\end{bmatrix}
\!+\!
\begin{bmatrix}0\\
0\\
\vdots\\
0\\
u_\mathrm{ref}
\end{bmatrix}
\label{eq:prod_laplace_state_space}
\end{equation}
The serial consensus form has several advantages, which will be the focus of the paper. First, however, we show that it satisfies the constraints we impose on the controller, as given by Definition~\ref{def:q_implementable}. In other words, we will discuss how the closed-loop structure in~\eqref{eq:prod_laplac_def} can be implemented on a network. 

When analysing the serial consensus controller of \eqref{eq:serial_consensus_controller} we will make use of the following assumption on the graph structure.
\begin{ass}(Connected spanning tree)\label{ass:connected_tree}
All graphs underlying the graph Laplacians $L_k$ contain a connected spanning tree.   
\end{ass}

\subsection{Implementing Serial consensus}
The following proposition ensures that the serial consensus system can be achieved by controlling the $n$\ts{th} order integrator system \eqref{eq:n_order_integrator} with an $n$-step implementable relative feedback controller as defined in Definition~\ref{def:q_implementable}.
\begin{proposition}\label{prop:q_implementable}
    Consider the $n$\ts{th}-order serial consensus as defined in \eqref{eq:prod_laplac_def}.
    If each $L_k\in \mathcal{A}^1(W,c)$ for some constant $c$ and adjacency matrix $W$, then the controller in \eqref{eq:serial_consensus_controller} is an %can be implemented with 
    $n$-step implementable relative feedback controller
    with respect to $W$ and a finite gain $c'$.
\end{proposition}
To prove this proposition we first need the following two lemmas whose proofs are provided in the appendix. 
\begin{restatable}{lemma}{sumlemma}\label{lem:q_sum}

    If $A_1\in\mathcal{A}^{q_1}(W,c_1)$ and $A_2\in\mathcal{A}^{q_2}~(W,c_2)$ then the sum $(A_1+A_2)\in\mathcal{A}^{\max(q_1,q_2)}~(W,c_1+c_2)$

\end{restatable}
\begin{restatable}{lemma}{prodlemma}\label{lem:q_product}

    Let $A_1\in\mathcal{A}^{q_1}(W,c_1)$ and $A_2\in \mathcal{A}^{q_2}(W,c_2)$ then the product $(A_1A_2)\in\mathcal{A}^{q_1+q_2}(W,c_1c_2)$

\end{restatable}
Now we can prove \propref{prop:q_implementable}.

\begin{proof}
    The serial consensus controller can be expanded to  the matrix polynomial%a sum
    \begin{align*}
    U(s)&=U_\mathrm{ref}(s)+\left(s^nI-\prod_{k=1}^{n} (sI+L_k)\right)X(s)\\
     &=U_\mathrm{ref}(s)+\left((s^n-s^n)I- \sum_{k=0}^{n-1} s^k A_k\right)X(s),     
    \end{align*}
    for some matrices $A_k$. To show the proposition, we need to show that $A_k \in \mathcal{A}^q(W,c')$ for all $k = 0,\ldots,n-1$, with $q\leq n$ and $c'<\infty$.
    Let 
    \begin{multline*}
        \mathcal{I}_k=\left\{\alpha~ \big| ~|\alpha|=n-k,~ \alpha \subset \{1,2,\dots,n \},\right.\\ 
    \left.i<j\implies \alpha(i)<\alpha(j) \right\}
    \end{multline*}
    denote all the ordered subsets of the range $[1,n]$ with size $n-k$. Then
    $$A_k=\sum_{\alpha\in \mathcal{I}_k}\prod_{j\in \alpha} L_j,\; \text{for all }k\in[0,n-1]. $$
    Since all $\alpha \in \mathcal{I}_k$ has $n-k$ elements we can show that $\prod_{j\in \alpha}L_j=B_\alpha\in \mathcal{A}^{n-k}~(W,c^{n-k})$ by applying Lemma~\ref{lem:q_product} recursively. Now we have a sum
    $$A_k=\sum_{\alpha\in \mathcal{I}_k} B_\alpha$$
    The number of ordered subsets of the range $[1,n]$ with size $n-k$ is given by the binomial coefficients and therefore the size of $|\mathcal{I}_k|=\binom{n}{n-k}$. Applying Lemma~\ref{lem:q_sum} recursively shows that $A_k\in \mathcal{A}^{n-k}(W,\binom{n}{n-k}c^{n-k})$. Clearly, we have that $n-k\leq n$ and $\binom{n}{n-k}c^{n-k}\leq \binom{n}{\lceil n/2 \rceil}\mathrm{max}(c,c^n)<\infty$ for all $k$. Let $c'=\binom{n}{\lceil n/2 \rceil}\mathrm{max}(c,c^n)$ and then we have that $A_k\in \mathcal{A}^n(W,c')$ for all $k$.    
\end{proof}
\begin{example}
    For clarity let us consider the controller for the case $n=3$. Then the controller is 
    \begin{multline*}        
        U(s)= U_\mathrm{ref}(s)+\left(s^3I-\prod_{k=1}^3 (sI+L_k)\right)X(s)\\
        =U_\mathrm{ref}(s)- \left(s^2(L_1+L_2+L_3)+\right. \\
        \left. s(L_1L_2+L_1L_3+L_2L_3)+L_1L_2L_3\right)X(s)
    \end{multline*}
    Here, $A_0=L_1L_2L_3$, $A_1=L_1L_2+L_1L_3+L_2L_3$, and $A_2=L_1L_2L_3$. The proposition asserts that if $L_1$, $L_2$, and $L_3$ share a sparsity pattern and have bounded gains, then the resulting controller gains $A_0$, $A_1$, and $A_2$ will be sparse and have bounded gains.
\end{example}

\subsection{Scalable stability}
Coordinating a multi-agent system is inherently a decentralized problem where the goal for each agent is to coordinate with its nearest neighbors. However when the controllers only depend on local measurements there is a possibility that controllers that manage to coordinate $N$ agents stop stabilizing as the number of agents increases. In~\cite{Tegling2022_Scalefrag} it was shown that for the $3$\ts{rd} and higher order consensus problem with controller $A_k=a_kL_N$ in~\eqref{eq:local_control},  the closed loop system will become unstable if the algebraic connectivity $\lambda_2(L_N)\rightarrow 0$ as $N\rightarrow \infty$. This motivates the notion of \emph{Scalable stability}
\begin{definition}[Scalable stability {\cite[Def.~2.1]{Tegling2022_Scalefrag}}]
A consensus control design is scalably stable if the resulting closed-loop system
achieves consensus over any graph in the family~$\{ \mathcal{G}_N \}$.% of finite size $N$.
\end{definition}

\section{Main Results}\label{sec:Results}
Our main contribution is two-fold. First we show the serial consensus achieves scalable stability and then we show that the implementation is robust to two classes of perturbations

\subsection{Scalable stability}

\begin{theorem}    
\label{thrm:high_order}
    Consider the $n$\ts{th} order serial consensus system as defined in \textit{Definition}~\ref{def:serial_consensus} under  Assumption~\ref{ass:connected_tree} and with $U_\mathrm{ref}\in \mathcal{R}\mathcal{H}_\infty$. Then the closed loop dynamics have the following properties:
    \begin{enumerate}[label=(\roman*)]
        \item The poles of \eqref{eq:prod_laplac_def} are given by the union of the eigenvalues of $-L_k$.  \label{prop:high_i}
        \item The solution achieves $n\ts{th}$ order consensus.
    \end{enumerate}
\end{theorem}
\vspace{1mm}
\begin{proof}
    \ref{prop:high_i} Any square matrix can be unitarily transformed to upper triangular form by the Schur traingularization theorem. Let $U_k L_kU_k^H= T_k$ be upper triangular. Then the block diagonal matrix $U=\mathrm{diag}(U_1, U_2, \dots U_{n})$ is a unitary matrix that upper triangularizes $A$ in \eqref{eq:prod_laplace_state_space}. For any triangular matrix the eigenvalues lie on the diagonal and this will be the eigenvalues of each $-L_k$.
    
    (ii) First, %let's 
    consider the closed loop dynamics of \eqref{eq:prod_laplac_def} which will be 
    $$X(s)=\left(\prod_{k=n}^1 (sI+L_{k})^{-1}\right)U_\mathrm{ref}(s).$$ Since, $U_\mathrm{ref}$ is stable, we know that the limit $\lim_{s\rightarrow 0} U_\mathrm{ref}(s)=U_\mathrm{ref}(0)$ exists. To prove that the system achieves $n\ts{th}$ order consensus we want to show that 
    $$\lim_{t\rightarrow \infty}y(t)=\lim_{s\rightarrow 0}C(s)X(s)=0$$
    for some transfer matrix $C(s)$, which encodes the consensus states. But since the reference dependence is only related to $U_\mathrm{ref}(0)$, we can simplify the problem to only consider impulse responses. But the impulse response has the same transfer function as the initial value response where $\xi_n(0)=U_\mathrm{ref}(0)$. Therefore, WLOG, assume that $U_\mathrm{ref}(s)=0$ and an arbitrary initial condition 
    $$\mathbf{\xi}(0)=[\xi_1(0)^T,\xi_2^T,\dots,\xi_n(0)^T]^T.$$ 
    The solution of~\eqref{eq:prod_laplace_state_space} is given by $\exp(At)\mathbf{\xi}(0)=S\exp(J(A)t)S^{-1}\mathbf{\xi}(0)$ where $J(A)$ is the Jordan normal form of $A$ and $S$ is an invertible matrix. From \ref{prop:high_i} and the diagonal dominance of the graph Laplacians we know that all eigenvalues of $A$ lie in the left half plane. By  Assumption~\ref{ass:connected_tree} it follows that the zero eigenvalue for each $L_k$ is simple. Now we prove that these $n$ zero eigenvalues form a Jordan block of size $n$. Let $e_k$ denote the $k$\ts{th} $1$-block vector, e.g. $\mathbf{e}_1=\begin{bmatrix}
    \mathbf{1}^T_N& 0_N & \dots & 0_N
    \end{bmatrix}^T$ and $e_2=\bmat{0_N& \mathbf{1}^T_N & 0_N&  \dots & 0_N}^T$. Then $e_1$ is an eigenvector since $A\mathbf{e}_1=0$. For $k\in\{2,3\dots,n\}$ we have $Ae_k=e_{k-1}$ which implies that $A^ke_k=0$. This shows that there is a Jordan block of size $n$ with an invariant subspace spanned by the vectors $e_k$. Since all other eigenvectors make up an asymptotically stable invariant subspace, it follows that $\xi(t)$ will converge towards a solution in $\mathrm{span}(e_1,e_2\ldots,e_n)$ and thus $\lim_{t\rightarrow \infty} \xi_k(t) =\alpha_k(t) \mathbf{1}_N$. Now, since $x(t)=\xi_1(t)$, it follows that $\lim_{t\rightarrow \infty}x(t)=\alpha_1(t)\mathbf{1}_N$, and furthermore, since 
    $$\dot{\xi}_k= -L_k\xi_k +\xi_{k+1}\rightarrow \xi_{k+1}\text{ as } t\rightarrow \infty$$
    for $k \in \{1,\dots,n-1\}$, it follows that $\lim_{t\rightarrow \infty}x^{(k)}(t)= \alpha_{k+1}(t) 1_N$ which shows that the system achieves $n$\ts{th} order consensus.
\end{proof}
\vspace{1mm}
This proposition shows that the stability of the consensus for the $n$\ts{th} order serial consensus can be reduced to verifying that the $n$ first-order consensus systems $\dot{x}=-L_k x$ achieve consensus.  %which 
This is equivalent to determining whether the graphs underlying each $L_k$ contains a connected spanning tree. This result together with \propref{prop:q_implementable} shows that $n$\ts{th} order consensus can be achieved with a local relative feedback controller with finite gain and thus achieve scalable stability. This result can be compared with \cite{Tegling2022_Scalefrag} where it is shown that the conventional consensus is not scalably stable for any order larger than $n=3$ if a graph Laplacian with vanishing algebraic connectivity is used. We can summarize this fact in the following corollary:
\begin{corollary} For any $n$, the controller~\eqref{eq:serial_consensus_controller} is scalably stable over any graph family~$\{\mathcal{G}_N\} $ that underlies $L_k$, provided  each $\mathcal{G}_N$ satisfies Assumption~\ref{ass:connected_tree}. 
\end{corollary}
\begin{rem}
Note that, by Theorem~4, scalable stability is achieved also with different graph families underlying each~$L_k$, and $||L_k||_\infty$ are allowed to be arbitrarily small. 
\end{rem}
\subsection{Robustness of serial consensus}
The proposed controller in \eqref{eq:serial_consensus_controller} is a relative state-feedback controller which is designed to ensure that the closed loop system achieves $n$\ts{th} order consensus as %which is 
guaranteed through \textit{Theorem}~\ref{thrm:high_order}. However, the $n$\ts{th} order integrator system may be an idealization of the system and the relative state feedback may need observers to be fully realized, and there could be unmodeled dynamics. These potential sources of errors call for a robust controller. We will now present two theorems, which prove that the serial consensus is robust towards two different types of uncertainties.
\subsubsection{Additive perturbation}
The following theorem asserts that the $n$\ts{th} order serial consensus controller can handle additive uncertainties.\begin{theorem}\label{thr:serial_robustness}
Consider the $n$\ts{th} order serial consensus system as defined in \textit{Definition}~\ref{def:serial_consensus}, under  Assumption~\ref{ass:connected_tree}, with $L_k=L$ for all $k$, and $L=L^T$. Then the perturbed system 
    $$(sI+L)^n X=U_\mathrm{ref}+\left(\sum_{k=0}^n\Delta_ks^kL^{n-k} \right)X,$$
    where $U_\mathrm{ref}, \Delta_k\in \mathcal{R}\mathcal{H}_\infty$, achieves $n$\ts{th} order consensus if  $$\|\Delta_0\|_{\mathcal{H}_\infty}+\|\Delta_n\|_{\mathcal{H}_\infty}+\sum_{k=1}^{n-1} \|\Delta_k\|_{\mathcal{H}_\infty}\sqrt{\dfrac{k^{k}}{n^{n}}(n-k)^{n-k}} <1$$
\end{theorem}
\vspace{2mm}

\begin{figure}
    \centering
    \includegraphics[width=0.9\linewidth]{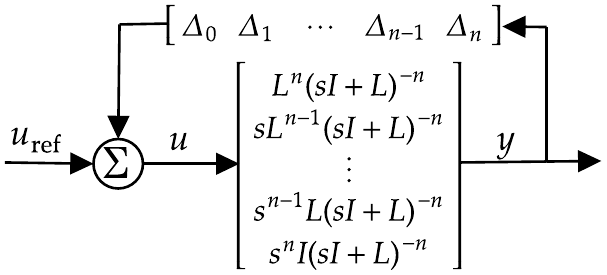}
    \caption{Block diagram illustrating the perturbation model in proof of \textit{Theorem} \ref{thr:serial_robustness}. }
    \label{fig:block_diagram_big}
\end{figure}

\begin{figure}
    \centering    \includegraphics[width=0.8\linewidth]{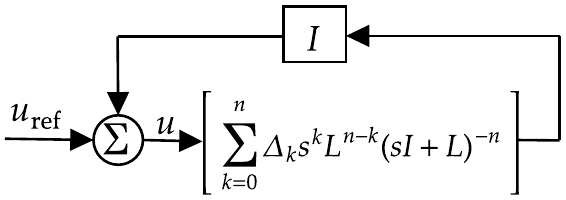}
    \caption{Block diagram illustrating the perturbation model in proof of \textit{Theorem} \ref{thr:serial_robustness}. }
    \label{fig:block_diagram_simpler}
\end{figure}

\begin{proof}
    First, note that the closed-loop system can be represented by the block diagram in \figref{fig:block_diagram_big}, which in turn can be simplified to \figref{fig:block_diagram_simpler}. Since $U_\mathrm{ref}$ is stable we can apply the small-gain theorem which asserts that $U(s)$ (as defined in the figures) will be stable if  
    $$\|\sum_{k=0}^n \Delta_ks^kL^{n-k}(sI+L)^{-n}\|_{\mathcal{H}_\infty}<1.$$
    Applying the triangle inequality and submultiplicativity on the left-hand side $(LH)$ yields
    \begin{equation}
        LH\leq \sum_{k=0}^n\|\Delta_k\|_{\mathcal{H}_\infty}\| s^kL^{n-k}(sI+L)^{-n}\|_{\mathcal{H}_\infty}
        \label{eq:LH_upperbound}
    \end{equation}
    Since $L$ is symmetric, it is possible to unitarily diagonalize it. Let $U=U^H$ denote one such unitary matrix. Then $L= U\Lambda U^H$ where $\Lambda$ is a non-negative real diagonal matrix. 
    $$\| s^kL^{n-k}(sI+L)^{-n}\|_{\mathcal{H}_\infty}=\| s^k\Lambda^{n-k}(sI+\Lambda)^{-n}\|_{\mathcal{H}_\infty}.$$
    For a diagonal matrix the singular values are given by the absolute value of the diagonal. Let, $\lambda>0$ be an arbitrary positive constant. The maximum gain for each diagonal can then be calculated through 
    $$\max_\omega \abs{\dfrac{\omega^k\lambda^{n-k}}{(j\omega+ \lambda)^n}}=\sqrt{\max_\omega \dfrac{\omega^{2k}\lambda^{2n-2k}}{(\omega^2+ \lambda^2)^n}}.$$
    The latter optimization problem is given by a continuous function and thus the derivative must be $0$ at the maximum. Simple calculus shows that the optimum is found at $\omega^2=\lambda^2k/(n-k)$ for $k=0, 1,\dots n-1$ and at $\omega=\infty$ for $k=n$. Inserting yields
    $$\max_\omega \abs{\dfrac{\omega^k\lambda^{n-k}}{(j\omega+ \lambda)^n}}=
    \left\{ \begin{matrix}
        \sqrt{\frac{k^k}{n^n} (n-k)^{n-k}} & \mathrm{if }\; 0<k<n\\
        1 & \mathrm{else } 
            \end{matrix}\right.        
    $$
    Now for the case where $\lambda=0$. Then we have for $k=0,\ldots, n-1$ 
    $$ \max_\omega \abs{\dfrac{\omega^k 0^{n-k}}{(j\omega+0)^n}}=0$$
    and for $k=n$
    $$\max_\omega \abs{\dfrac{\omega^n}{(j\omega+0)^n}}=1.$$
    This is less restrictive than for $\lambda>0$ and thus we can use the result for $\lambda>0$. Plugging this into the upper bound of the $LH$ \eqref{eq:LH_upperbound} results in the sought inequality. 

    Finally, we must ensure that stability of the closed loop in \figref{fig:block_diagram_simpler} implies $n$\ts{th} order consensus. Since the transfer matrix from $u$ to $y$ in \figref{fig:block_diagram_big} is stable it follows that $Y(s)$ will be stable. This means that we have shown the following $\lim_{t\rightarrow\infty} L^{n-k}x^{(k)}(t)= 0$. By Assumption~\ref{ass:connected_tree} the $0$ eigenvalue of $L$ is unique and therefore $0$ is a unique eigenvalue of $L^{n-k}$ too. Subsequently, $\lim_{t\rightarrow \infty} x^{(k)}(t)\in \mathrm{ker}(L^{n-k})$. Since $L^{n-k}\mathbf{1}_N=0$ it follows that $\lim_{t\rightarrow \infty} x^{(k)}(t)\in\mathrm{span}(\mathbf{1}_N)$ and that the agents will reach consensus in all the $n-1$ first time derivatives and thus achieve $n$\ts{th} order consensus.
\end{proof}
It is worth noting that the norm bound on the uncertainty blocks $\Delta$ is independent of the number of agents in the system. Therefore, the serial consensus implementation is \emph{scalably robust} in the sense that it allows equally sized perturbations regardless of network size. This is not the case for localized conventional consensus, following the results in~\cite{Tegling2022_Scalefrag}.

\subsubsection{ Multiplicative perturbation}
It is also possible to see the closed-loop serial consensus system as a series of interconnected first-order systems. Therefore it is also interesting to consider the robustness with respect to the individual factors. The following theorem gives a sufficient condition  for the unforced closed loop system to achieve $n$\ts{th} order consensus.
\begin{theorem}\label{thr:robustness_series}
    The following perturbed $n$\ts{th} order serial consensus system
    $$
        \left( sI+s\Delta_0+(I+\Delta_n)L_n\right)
        \prod_{k=1}^{n-1}\left(sI+(I+\Delta_k)L_k\right) X=U_\mathrm{ref}$$
        where $U_\mathrm{ref},\Delta_k\in \mathcal{R}\mathcal{H}_\infty$ and $L_k=L_k^T$ for $k = 1,\ldots,n$,
    achieves $n$\ts{th} order consensus if 
    $$\|\Delta_k\|_{\mathcal{H}_\infty}<1, \text{ for all } k$$
    and
    $$\|\Delta_0\|_{\mathcal{H}_\infty}+\|\Delta_n\|_{\mathcal{H}_\infty}<1.$$
 
\end{theorem}
\vspace{2mm}

\begin{figure}
    \centering
    \includegraphics[width=0.9\linewidth]{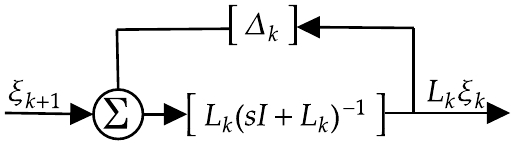}
    \caption{Block diagram illustrating the perturbation model of a general first-order consensus block which is used in the proof of \textit{Theorem} \ref{thr:robustness_series}.}
    \label{fig:block_diagram_small}
\end{figure}

\begin{proof}
    First, note that we can construct $X(s)=\Xi_1(s)$ and $s{\Xi}_k=-(I+\Delta_k)L_k\Xi_k+\Xi_{k+1}$ for $k=1,\dots,n-1$ and $s(I+\Delta_0)\Xi_n=-(I+\Delta_n)L_n\Xi_n+ U_\mathrm{ref}$. For $\Xi_n$ we have exactly the $1$\ts{st} order case of \textit{Theorem}~\ref{thr:serial_robustness} and thus $\lim_{t\rightarrow\infty}\xi_n(t)=\alpha_n(t)\mathbf{1}_N$ if  $\|\Delta_0\|_{\mathcal{H}_\infty}+\|\Delta_n\|_{\mathcal{H}_\infty}<1$. Consider the following induction hypothesis: if $\Xi_{k+1}(s)=\mathbf{1}_NG_{k+1}(s)+H_{k+1}(s)$ where $H_{k+1}(s)\in \mathcal{R}\mathcal{H}_\infty$ then $\Xi_{k}=\mathbf{1}_NG_{k}(s)+H_{k}(s)$  for some $H_{k}(s)\in \mathcal{R}\mathcal{H}_\infty$. We have 
    $$s{\Xi}_k=-(I+\Delta_k)L_k\Xi_k+\Xi_{k+1}$$
    which can be represented by the block diagram \figref{fig:block_diagram_small}. Here, note that $$L_k(sI+L_k)^{-1}\Xi_{k+1}=(sI+L_k)^{-1}L_k(H_{k+1}(s))$$
    and the potentially unstable term of $\Xi_{k+1}$ can be ignored. Reusing a result from the previous proof we have $\|L_k(sI+L_k)^{-1}\|_{\mathcal{H}_\infty}=1$ and therefore $L_k\Xi_k \in \mathcal{R}\mathcal{H}_\infty$ if $\|\Delta_k\|_{\mathcal{H}_\infty}<1$. Since the $0$ eigenvalue of $L_k$ is unique, it follows that $\Xi_k(s)=\mathbf{1}_NG_{k}(s)+H_{k}(s)$ with $H_k\in\mathcal{R}\mathcal{H}_\infty$ which proves the induction hypothesis since we have already shown the base case $\Xi_n(s)=\mathbf{1}_NG_n(s) +H_n(s)$.
    Left is to prove that the system will reach $n$\ts{th} order consensus. Note that $L_1X(s)=L_1\Xi_1(s)$ is stable and therefore we get through the final value theorem
    $$\lim_{t\rightarrow \infty}L_1x(t)=\lim_{s\rightarrow 0}sL_1\Xi_1(s)=0.$$
    Furthermore, we have for all $k$:    $\lim_{s\rightarrow 0}sL_k\Xi_k(s)=0$.
    This, combined with $s^2\Xi_k(s)=-(I+\Delta_k)sL_k\Xi_k(s)+s\Xi_{k+1}(s)$
    shows that 
\begin{align*}
\lim_{t\rightarrow \infty} L_{k+1}x^{(k)}(t)=\lim_{s\rightarrow 0} s (s^{k} L_{k+1}X(s))\\
    =\lim_{s\rightarrow 0} sL_{k+1}\Xi_{k+1}(s)=0
\end{align*}
Finally, since the $0$ eigenvalues for each $L_k$ are unique with corresponding eigenvector $\mathbf{1}_N$ we see that $n$\ts{th} order consensus will be achieved.
\end{proof}
This theorem shows that the $n$\ts{th} order serial consensus is robust in its construction.

\section{Examples}\label{sec:Examples}
\subsection{$2$\ts{nd} order consensus on circular graph}
Consider the directed cycle graph, which is represented by the adjacecncy matrix 
$$[W]_{i,j}=1 \text{ iff } i-j=1 \mod{N}.$$
The corresponding graph Laplacian $L_c$ is a circulant matrix and therefore, the eigenvalues are known analytically. In particular, the eigenvalue with the second smallest real part is
$\lambda_2(L_c)=1-\exp(2\pi/N)=1-\cos{(2\pi/N)}-\mathbf{i}\sin{(2\pi/N)}$. For large $N$, this eigenvalue can be approximated with a $1$\ts{st} order Taylor approximation, which yields
$\lambda_2(L_c)\approx -\mathbf{i} 2\pi/N$. This eigenvalue will cause problems when designing a controller using the conventional consensus. To see this, consider the closed loop dynamics
$$s^2I+2p_1sL_c+p_0L_c=U_\mathrm{ref}.$$
The system can be diagonalized and, in particular, two of the poles are given by the equation 
$$s^2+2p_1\lambda_2(L_c)+p_0\lambda(L_c) =0.$$
In the case that $p_0$ and $p_1$ are designed independently of the network size $N$, then for sufficiently large $N$ the roots %solution 
can be approximated as
$$s_p=-p_1 \lambda_2 \pm \sqrt{p_1^2\lambda_2^2-p_0\lambda}\approx\pm (1+\mathbf{i})\sqrt{\frac{\pi p_0}{N}}  $$
Since one of these poles will eventually lie in the RHP, it follows that the closed loop system will become unstable when $N$ is sufficiently large, regardless of the choice of $p_0$ and $p_1$.

For the serial consensus it suffices to check that all eigenvalues but the unique $0$ eigenvalue of $L_c$ lie in the RHP or equivalently if $0<\mathrm{Re}(\lambda_2(L_C))=1-\cos(2\pi/N)$ which is clearly true for any finite $N$. Alternatively, it is also sufficient that the underlying graph contains a connected spanning tree.

\subsection{3\ts{rd} order consensus}
\begin{figure}
     \centering
    \begin{subfigure}[]{.48\linewidth}
        \centering
      \includegraphics[width=\linewidth]{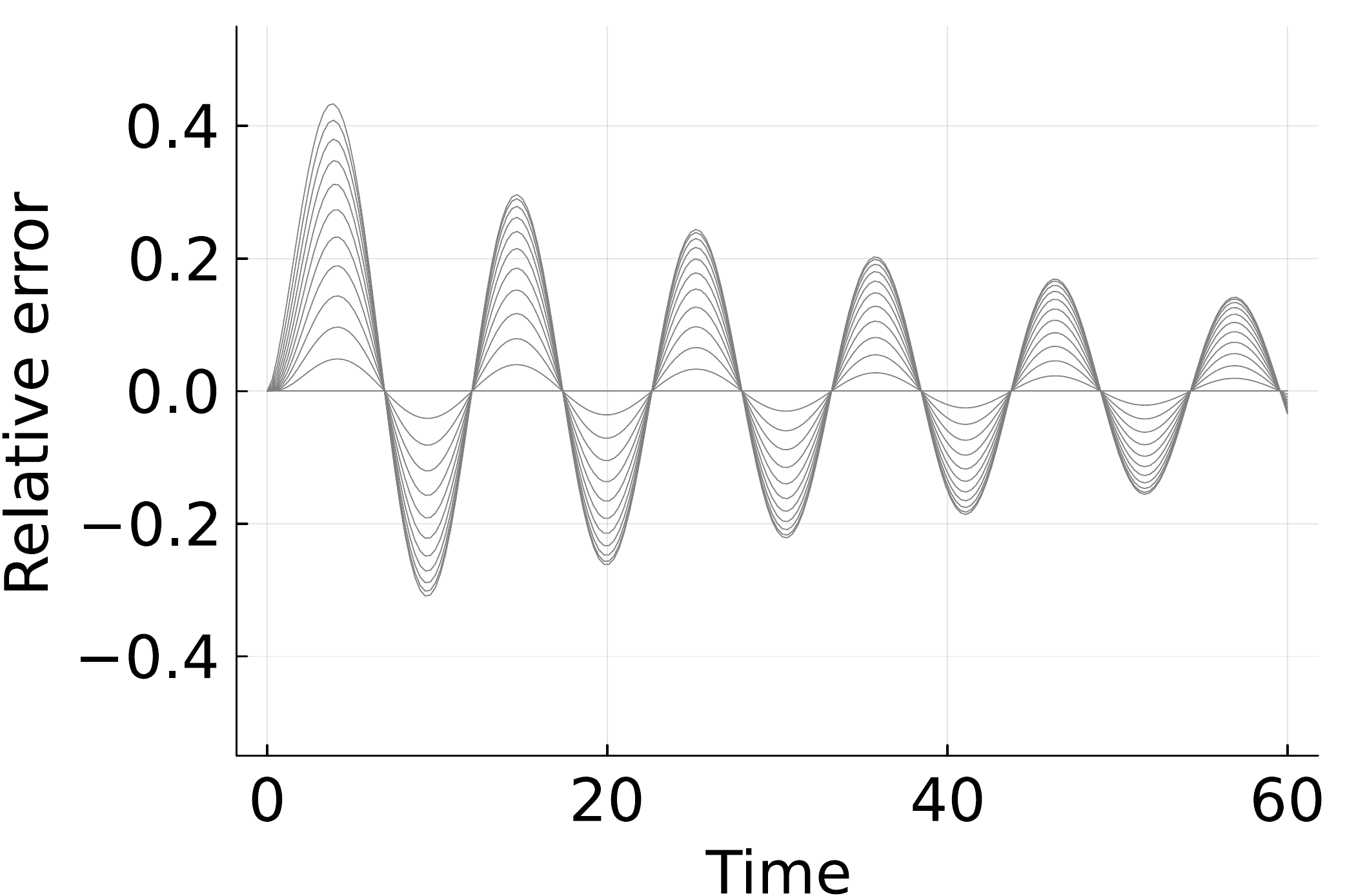}   
        \caption{Conventional with $N=12$.}
        \label{fig:conv12}
    \end{subfigure}
        \hfill
    \begin{subfigure}[]{.48\linewidth}
        \centering
        \includegraphics[width=\linewidth]{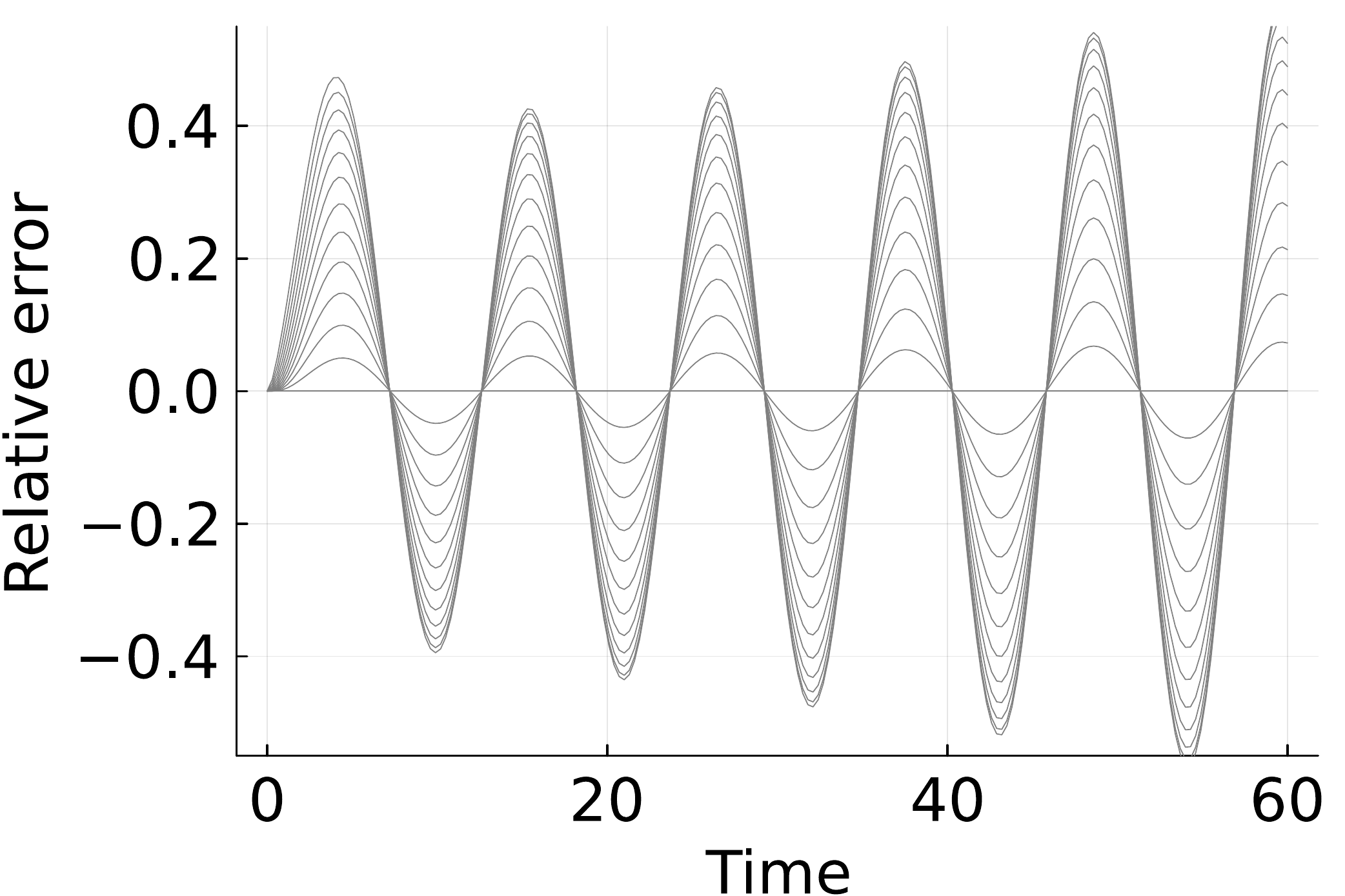}
        \caption{Conventional with $N=13$.}
        \label{fig:conv13}
    \end{subfigure}
    \medskip
    
    \begin{subfigure}[]{.48\linewidth}
    \centering
    \includegraphics[width=\linewidth]{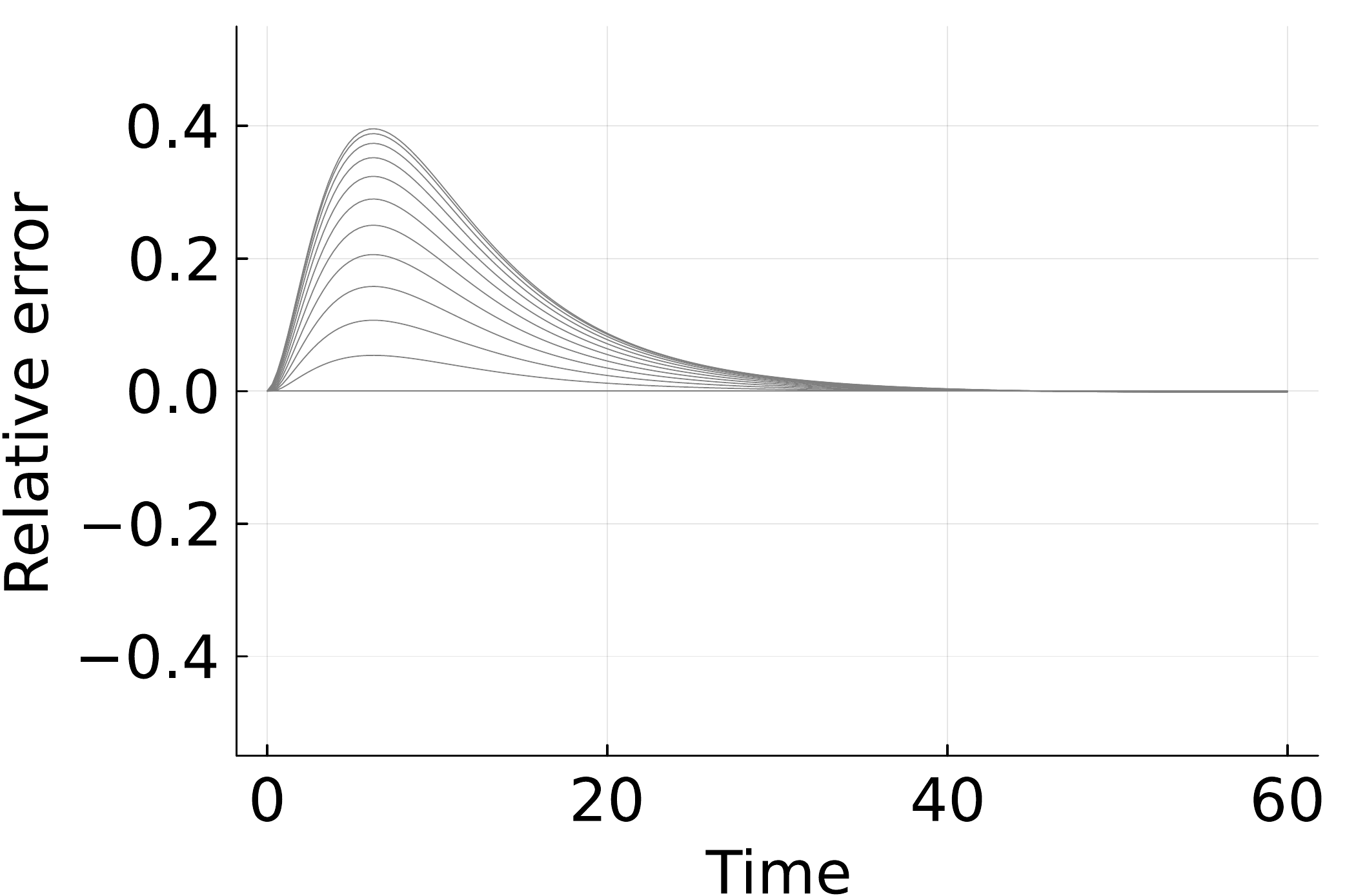}
    \caption{Serial with $N=12$.} 
    \label{fig:ser12}
    \end{subfigure}
       \hfill
    \begin{subfigure}[]{.48\linewidth}
        \centering
        \includegraphics[width=\linewidth]{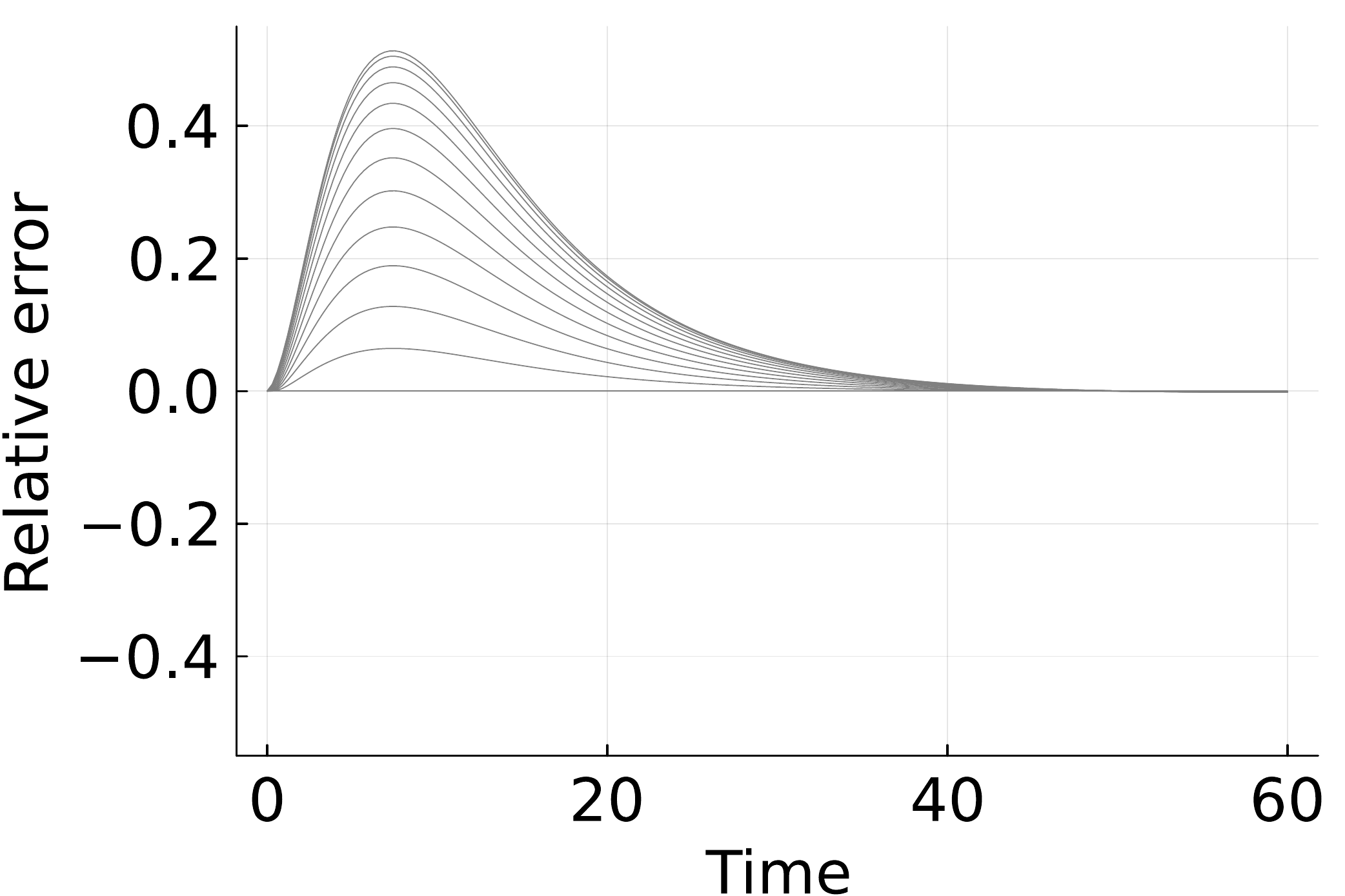}
        \caption{Serial with $N=13$}
        \label{fig:ser13}
    \end{subfigure}
     \caption{$3$\ts{rd} order consensus in a chain of vehicles is considered. The plots show the intervehicle relative errors over time when the lead vehicle moves at constant acceleration. Panels (a) and (b) show that the addition of one agent destabilizes the closed loop for the conventional consensus. Panels (c) and (d) illustrate the fact that the serial consensus will remain stable under such agent additions.}
     \label{fig:Conv_serial}
\end{figure}

It has been shown that for $n\geq 3$ it is not possible to achieve scalable stability for any graph family $\{G_N\}$ where the corresponding graph Laplacian $L_N$ has an eigenvalue with vanishing real part as the graph is growing, i.e. 
 if $\lim_{N\rightarrow \infty}\mathrm{Re}(\lambda_2(L_N)) =0$. At least, this is not possible with the conventional consensus control. On the serial consensus form this is no longer a problem. The controller 
 $$U(s)= U_\mathrm{ref}+\left(s^3I-\prod_{k=1}^{3} (sI+L_N) \right)X(s)$$
will achieve consensus as long as the underlying graphs~$\{G_N\}$ all contains a connected spanning tree. To illustrate this, consider the graph defined by $W\in \mathbb R^{N\times N}$, the adjacency matrix
$$W_{i,j}=\left\{ \begin{matrix} 1 & \text{if }|i-j|=1\text{ and } i\neq 1\\
0 & \text{else}    
\end{matrix}\right. .$$
This corresponds to a bidirectional string with a leader (Agent~$1$). Let $L$ be the associated graph Laplacian. 
It is true that $\lim_{N\rightarrow \infty} \lambda_2(L)= 0$ and thus any conventional control design with $L$ will eventually lead to an unstable closed loop. For this example, let the conventional control law be
$u(t)=u_\mathrm{ref}(t)-6L\ddot{x}-4L\dot{x}-2Lx$
and the serial consensus controller \eqref{eq:serial_consensus_controller} be defined with the same graph Laplacians $L_k=2kL$. The response to a constant acceleration of the leader is shown in Fig.~\ref{fig:Conv_serial}. Here we see that the addition of a $13$th agent to the system destabilizes the closed loop for the conventional consensus while the serial consensus only loses some performance.

\subsection{Robustness of the $2$\ts{nd} order serial Consensus.}
Theorems~\ref{thr:serial_robustness}~and~\ref{thr:robustness_series} show that the serial consensus can be perturbed and still achieve $n$\ts{th} order consensus. Now we want to illustrate what the block $\Delta_k$ can be. Consider the perturbed $2$\ts{nd} order consensus system in Theorem~\ref{thr:serial_robustness}. Writing out all terms we get
$$s^2(I+\Delta_2)X=U_\mathrm{ref}-(s (2I+\Delta_1)LX+(I+\Delta_0)L^2X).$$
In this form, the $\Delta_2$ block can be thought of as representing model errors; we may control a system which we model as being $N$ identical double integrator systems but in reality they may differ. This is obviously the case for vehicle platoons, which are often modeled as chains of identical double integrators. Through our theorem we can for instance allow $\Delta_2$ to be a diagonal transfer matrix with elements $[\Delta_2]_{i,i}=\frac{k_i}{T_is+1}$ where $|k_i|<1$ and $T_i>0$ for all $i$. Then, the closed loop system would remain stable despite the heterogeneous agents. The blocks $\Delta_1$ and $\Delta_0$ are also important. For instance the signals $L^2x(t)$ and $L\dot{x}(t)$ may not be directly measured but estimated through linear filters. This could be thought of as unmodeled dynamics which these blocks can capture.

If we focus on Theorem~\ref{thr:robustness_series}, then the perturbed model is
$$(s(\Delta_0+I)+(\Delta_1+I)L_1)(sI+(\Delta_2+I)L_2)X=U_\mathrm{ref}$$
The theorem only asserts robustness for symmetrical graph Laplacians $L_k$. However, since $\Delta_k$ can also be constant matrices, it is also possible to construct new  (asymmetric) graph Laplacians $L_k'=(I+\Delta_k)L_k$ by designing the $\Delta_k$ blocks. 

\section{Conclusion}\label{Conclusion}

This work has introduced the $n$\ts{th} order serial consensus system which can be seen as a natural generalization of the well-known consensus protocols. The stability of the introduced system can be analysed by considering $n$ regular first order consensus protocol. The proposed controller to achieve $n$\ts{th} order serial consensus has been shown to be implementable using relative measurements confined to a local neighborhood of each agent and can therefore be considered a decentralized control scheme. Robustness of the proposed system has also been analyzed. This has been addressed in terms of additive and model perturbations. The analysis showed that the size, measured in the $\mathcal{H}_\infty$ norm, of the allowable uncertainties were independent of the number of agents.

Future and ongoing work includes looking into the performance  of the serial consensus and how this relates to string stability. It would also be interesting to look into an implementation where each agent implements an observer to compute their control action.

\appendix \label{app:lemma_proofs}

Here we prove Lemmas~\ref{lem:q_sum},~\ref{lem:q_product} which describe how the sparsity pattern of two matrices changes through addition and multiplication. The Lemmas are restated for convenience.
\sumlemma*
\begin{proof}
First, we have $\|A_1+A_2\|_\infty\leq \|A_1\|_\infty+\|A_2\|_\infty\leq c_1+c_2$ which follows from the triangle inequality.

For the second part we have $(A_1+A_2) \mathbf{1}_N=0 +0=0$. 

For the last part, WLOG, suppose that $q_1\leq q_2=\max(q_1,q_2)$. Since $W$ is a positive matrix, we get
$$0\leq \left(\sum_{k=0}^{q_1}W^k\right)_{i,j}\leq\left(\sum_{k=0}^{q_2}W^k\right)_{i,j}.$$
In particular, the following implication follows
$$ \left(\sum_{k=0}^{q_2}W^k\right)_{i,j}=0 \implies \left(\sum_{k=0}^{q_1}W^k\right)_{i,j}\implies [A_1+A_2]_{i,j}=0$$
\end{proof}

To prove the result on the product of two matrices, Lemma~\ref{lem:q_product}, we need the following three lemmas:
\begin{lemma}\label{lem:abs_change}
    Let $A,B\in\mathbb{C}^{N\times N}$ and define $\hat{A}_{i,j}=|A|_{i,j}$ and $\hat{B}_{i,j}=|B|_{i,j}$. If $(A B )_{i,j}\neq 0$ then $(\hat{A}\hat{B})_{i,j}\neq 0$.
\end{lemma}
\begin{proof}
    Suppose the statement is false, i.e. $(\hat{A}\hat{B})_{i,j}=0$ but $(AB)_{i,j}\neq0$. Then we know that 
    $$(\hat{A}\hat{B})_{i,j}=\sum_{k=1}^N |A_{i,k}||B_{k,j}|=0,$$
    but this implies that at least one of $A_{i,k}$ and $B_{k,j}$ is equal to $0$ for all $k$. But from this it follows that
    $$(AB)_{i,j}=\sum_{k=1}^N A_{i,k}B_{k,j}=\sum_{k=1}^N 0=0.$$
    This is a contradiction and concludes the proof.
\end{proof}

\begin{lemma}\label{lem:element_addition}
    Let $A,A_1,B,B_1\in\mathbb{R}_+^{N\times N}$. If $(AB)_{i,j}\neq 0$ then
    $\left( (A+A_1)(B+B_1)\right)_{i,j}\neq 0$. 
\end{lemma}
\begin{proof}
    Expand the product to get 
    \begin{multline*}        
    \left((A+A_1)(B+B_1)\right)_{i,j}=(AB)_{i,j}+(AB_1)_{i,j}\\+(A_1B)_{i,j}+(A_1B_1)_{i,j}\geq (AB)_{i,j}\end{multline*}
    which followed from the fact that the product of $2$ nonnegative matrices is also nonnegative.
\end{proof}

\begin{lemma}\label{lem:elem_rescale}
    Let $A,A_1,B,B_1\in\mathbb{R}_+^{N\times N}$ be such that $A_{i,j}=0$ if and only if ${A_1}_{i,j}=0$, and $B_{i,j}=0$ if and only if ${B_1}_{i,j}=0$. Then, $(AB)_{i,j}=0$ if and only if $(A_1B_1)_{i,j}=0$. 
\end{lemma}
\begin{proof}
    The statement is clearly symmetrical and it is enough to prove sufficiency. Now, if $(AB)_{i,j}=0$ then we know that
    $$\sum_k A_{i,k}B_{k,j}=0, \implies A_{i,k}B_{k,j}=0, ~\forall k$$
    But this implies that either $A_{i,k}=0$ or $B_{k,j}=0$. In turn, this implies that either $(A_1)_{i,k}=0$ or $(B_1)_{k,j}=0$. And this leads to
    $$(A_1B_1)_{i,j}=\sum_k (A_1)_{i,k}(B_1)_{k,j}=0$$
\end{proof}

Now we can prove Lemma~\ref{lem:q_product}:
\prodlemma*
\begin{proof}
First, the gain can be bounded as $\|A_1A_2\|_\infty \leq \|A_1\|_\infty\|A_2\|_\infty\leq c_1c_2$ which followed from submultiplicity of the induced norm and from the definition of the sets.

For the second part we have $A_1A_2 \mathbf{1}_N=A_10=0$. 

For the last part we have to do slightly more. First replace each element in $A_1$ and $A_2$ with its absolute value and denote these $B_1$ and $B_2$. Now introduce two non-negative matrices $C_1$ and
$C_2$ such that $B_1+C_1=0 \iff \sum_{k=0}^{q_1}W^k$ and $B_2+C_2=0 \iff \sum_{k=0}^{q_2}W^k$. Finally note that 
$$(\sum_k^{q_1}W^k)(\sum_j^{q_2} W^j)=\sum_k^{q_1+q_2}w_k W^k $$
for some $w_k>0$. By applying Lemma~\ref{lem:elem_rescale} two times we get that 
$$\left[\sum_k^{q_1+q_2}W^k\right]_{i,j}=0\implies \left[(B_1+C_1)(B_2+C_2)\right]_{i,j}=0$$
Through Lemma \ref{lem:element_addition} we get
$$\left[(B_1+C_1)(B_2+C_2)\right]_{i,j}=0 \implies \left[B_1B_2\right]_{i,j}=0$$
And finally applying Lemma \ref{lem:abs_change} results in
$$\left[B_1B_2\right]_{i,j}=0\implies \left[A_1A_2\right]_{i,j}=0$$
\end{proof}

\section*{Acknowledgement}
We want to thank Richard Pates for useful discussions regarding the robustness results.

\bibliographystyle{IEEETran}
\bibliography{references}   

\begin{thebibliography}{10}
\providecommand{\url}[1]{#1}
\csname url@rmstyle\endcsname
\providecommand{\newblock}{\relax}
\providecommand{\bibinfo}[2]{#2}
\providecommand\BIBentrySTDinterwordspacing{\spaceskip=0pt\relax}
\providecommand\BIBentryALTinterwordstretchfactor{4}
\providecommand\BIBentryALTinterwordspacing{\spaceskip=\fontdimen2\font plus
\BIBentryALTinterwordstretchfactor\fontdimen3\font minus
  \fontdimen4\font\relax}
\providecommand\BIBforeignlanguage[2]{{%
\expandafter\ifx\csname l@#1\endcsname\relax
\typeout{** WARNING: IEEEtran.bst: No hyphenation pattern has been}%
\typeout{** loaded for the language `#1'. Using the pattern for}%
\typeout{** the default language instead.}%
\else
\language=\csname l@#1\endcsname
\fi
#2}}

\bibitem{FaxMurray}
J.~A. Fax and R.~M. Murray, ``Information flow and cooperative control of
  vehicle formations,'' \emph{IEEE Trans. Autom. Control}, vol.~49, no.~9, pp.
  1465--1476, Sep 2004.

\bibitem{OlfatiSaber2004}
R.~Olfati-Saber and R.~M. Murray, ``Consensus problems in networks of agents
  with switching topology and time-delays,'' \emph{IEEE Trans. Autom. Control},
  vol.~49, no.~9, pp. 1520--1533, Sept 2004.

\bibitem{Jadbabaie2003}
A.~Jadbabaie, J.~Lin, and A.~S. Morse, ``Coordination of groups of mobile
  autonomous agents using nearest neighbor rules,'' \emph{IEEE Trans. Autom.
  Control}, vol.~48, no.~6, pp. 988--1001, June 2003.

\bibitem{Pasqualetti2014}
F.~Pasqualetti, S.~Zampieri, and F.~Bullo, ``Controllability metrics,
  limitations and algorithms for complex networks,'' in \emph{2014 American
  Control Conf.}, 2014, pp. 3287--3292.

\bibitem{Bamieh2012}
B.~Bamieh, M.~R. Jovanovi\'c, P.~Mitra, and S.~Patterson, ``Coherence in
  large-scale networks: {D}imension-dependent limitations of local feedback,''
  \emph{IEEE Trans. Autom. Control}, vol.~57, no.~9, pp. 2235--2249, Sep. 2012.

\bibitem{SiamiMotee2015}
M.~Siami and N.~Motee, ``Fundamental limits and tradeoffs on disturbance
  propagation in large-scale dynamical networks,'' \emph{IEEE Trans. Autom.
  Control}, vol.~61, no.~12, pp. 4055--4062, 2016.

\bibitem{swaroop1996stringstability}
D.~Swaroop and J.~Hedrick, ``String stability of interconnected systems,''
  \emph{IEEE Trans. Autom. Control}, vol.~41, no.~3, pp. 349--357, March 1996.

\bibitem{seiler2004disturbancep_propagation}
P.~Seiler, A.~Pant, and K.~Hedrick, ``Disturbance propagation in vehicle
  strings,'' \emph{IEEE Trans. Autom. Control}, vol.~49, no.~10, pp.
  1835--1842, Oct 2004.

\bibitem{Tegling2022_Scalefrag}
E.~Tegling, B.~Bamieh, and H.~Sandberg, ``Scale fragilities in localized
  consensus dynamics,'' \emph{Automatica}, 2023, to appear. Available:
  https://arxiv.org/abs/2203.11708.

\bibitem{Jiang2009}
F.~Jiang, L.~Wang, and Y.~Jia, ``Consensus in leaderless networks of
  high-order-integrator agents,'' \emph{American Control Conf.}, pp.
  4458--4463, June 2009.

\bibitem{Ren2006}
W.~Ren, K.~Moore, and Y.~Chen, ``High-order consensus algorithms in cooperative
  vehicle systems,'' \emph{IEEE International Conf. on Networking, Sensing and
  Control}, pp. 457--462, 2006.

\bibitem{Studli2017}
S.~Stüdli, M.~M. Seron, and R.~H. Middleton, ``Vehicular platoons in cyclic
  interconnections with constant inter-vehicle spacing,''
  \emph{IFAC-PapersOnLine}, vol.~50, no.~1, pp. 2511--2516, 2017, 20th IFAC
  World Congress.

\bibitem{zhou1998essentials}
K.~Zhou and J.~C. Doyle, \emph{Essentials of robust control}.\hskip 1em plus
  0.5em minus 0.4em\relax Prentice hall Upper Saddle River, NJ, 1998, vol. 104.

\end{thebibliography}
\vspace{1mm}
\end{document}